\documentclass{article}%
\usepackage{amssymb}
\usepackage{amsmath}
\usepackage{amsfonts}
\usepackage{graphicx}%
\newtheorem{theorem}{Theorem}

\newtheorem{corollary}[theorem]{Corollary}

\newtheorem{definition}[theorem]{Definition}

\newtheorem{lemma}[theorem]{Lemma}

\newtheorem{proposition}[theorem]{Proposition}
\newtheorem{remark}[theorem]{Remark}

\begin{document}

\title{Hylomorphic solitons for the Benjamin-Ono and the fractional KdV equations}
\author{Vieri Benci\thanks{ Dipartimento di Matematica, Universit\`{a} degli Studi di
Pisa, Via F. Buonarroti 1/c, Pisa, ITALY\texttt{ }and Centro Interdisciplinare
"Beniamino Segre", Accademia dei Lincei. e-mail: \texttt{vieri.benci@unipi.it}%
}, Donato Fortunato\thanks{Dipartimento di Matematica, Universit\`{a} di Bari
"Aldo Moro" Via Orobona 4, 70125 Bari and INFN, Sezione 4 - email:
\texttt{donato.fortunato@uniba.it}}}
\maketitle

\begin{abstract}
This paper concerns with the existence of solitons, namely stable solitary
waves, for the Benjamin-Ono and the fractional KdV equations.\bigskip

\textbf{AMS subject classification}: 74J35, 35C08, 35A15, 35Q74, 35B35

\bigskip

\textbf{Key words}: Benjamin-Ono equation, fractional KdV equation, fractional
Schroedinger equation, travelling solitary waves, hylomorphic solitons,
variational methods.\bigskip

\textit{Dedicated to the memory of Enrico Magenes.}

\end{abstract}
\tableofcontents

\section{Introduction}

The Benjamin-Ono equation (BO) is a model of one dimensional waves in deep
water
\begin{equation}
\partial_{t}u+\mathcal{H}\partial_{x}^{2}u+u\partial_{x}u=0 \tag{BO}\label{BO}%
\end{equation}
where $u=u(t,x),\ $and $\mathcal{H}$ denotes the spatial Hilbert transform:%
\[
\mathcal{H}w(x)=\frac{1}{\pi}p.v.\int\frac{w(y)}{x-y}dy=\frac{-i}{\sqrt{2\pi}%
}\int\frac{\xi}{|\xi|}\hat{w}(\xi)e^{ix\xi}d\xi
\]
where $\hat{w}(\xi)$ denotes the Fourier tranform of $w$, namely%
\[
\hat{w}(\xi)=\frac{1}{\sqrt{2\pi}}\int w(y)e^{-iy\xi}dy.
\]

It is well known that (\ref{BO}) admits soliton solutions (see e.g.
\cite{Kenig2009}). In this paper, we shall use the method developed in
\cite{befolib} to prove that a large class of equations including equation
(\ref{BO}) admits \textit{hylomorphic }solitons. Following \cite{milan} and
\cite{befolib}, a soliton is called \textit{hylomorphic} if its stability is
due to a particular interplay between the \textit{energy }$E$ and the
\textit{hylenic} \textit{charge }%
\[
C:=\int u^{2}dx
\]
which is another integral of motion. More precisely, a soliton $u_{0}$ is
hylomorphic if%
\[
E(u_{0})=\min\left\{  E(u)\ |\ \int u^{2}dx=C(u_{0})\right\}  .
\]
In this paper we give a general theorem which, if it is applied to (\ref{BO}),
gives the following theorem:

\begin{theorem}
\label{martina}Equation \ref{BO} has a one parameter family $u_{\delta
},\ \delta\in\left(  0,\delta_{\infty}\right)  ,$ of hylomorphic solitons (see
Def. \ref{tdc}). Moreover, $u_{\delta}$ is a (weak) solution of the equation%
\begin{equation}
\mathcal{H}\partial_{x}^{2}u+u\partial_{x}u=\lambda_{\delta}\partial_{x}u
\end{equation}
and
\begin{equation}
U_{\delta}(t,x)=u_{\delta}(x-\lambda_{\delta}t)
\end{equation}
solves (\ref{BO}) for suitable $\lambda_{\delta}$.
\end{theorem}

It is well known that, in this case (see e.g. \cite{Kenig2009}), eq.
(\ref{BO}) has explicit solutions, namely
\[
u_{\delta}(x-\lambda_{\delta}t)=\frac{4\lambda_{\delta}}{1+\lambda_{\delta
}^{2}(x-x_{0}-\lambda_{\delta}t)^{2}}.
\]

We get the above theorem as a particular case of the study of the following
two families of equations (see Th. \ref{marina} and Th. \ref{marina+}):%

\begin{equation}
\partial_{t}u+\partial_{x}\left[  D_{x}^{2s}u+W^{\prime}(u)\right]
=0,\ \ u(t,x)\in\mathbb{R},\ s\in\mathbb{R},\ s\geq\frac{1}{2} \tag{FKdV}%
\label{GBO}%
\end{equation}
and%
\begin{equation}
i\frac{\partial\psi}{\partial t}=\frac{1}{2}D^{2s}\psi+\frac{1}{2}W^{\prime
}(|\psi|)\frac{\psi}{|\psi|},\ \ \psi(t,x)\in\mathbb{C},\ s\in\mathbb{R}%
,\ s\geq\frac{1}{2} \tag{FNS}\label{GNS}%
\end{equation}
where
\begin{equation}
D^{s}w(x)=\frac{1}{\sqrt{2\pi}}\int|\xi|^{s}\hat{w}(\xi)e^{ix\xi}d\xi;
\label{df}%
\end{equation}
and $W\in C^{2}(\mathbb{R)}$. We will refer to these equations as to the
Fractional Korteweg--de Vries equation (FKdV) and the Fractional Nonlinear
Schroedinger (FNS) equation respectively. Moreover, it is immediate to see
that
\[
D=\mathcal{H}\partial_{x}.
\]
The above equations, for particular choises of $s$ and $W,$ reduce to well
known PDE's of physics. In particular if $W(r)=\frac{1}{6}r^{3}$ and $s=1/2$,
(\ref{GBO}) reduces to (\ref{BO}).

If $s=1,\ $and $W(r)=-\frac{1}{6}r^{3},\ $then$\ \ D^{2s}=D^{2}=-\partial
_{x}^{2}$ and hence (\ref{GBO}) reduces to the KdV equation:%

\begin{equation}
\partial_{t}u-\partial_{x}^{3}u-u\partial_{x}u=0 \label{1K}%
\end{equation}
When $s=m\in\mathbb{N},$ and $W(r)=-\frac{\left\vert r\right\vert ^{p}%
}{p\left(  p-1\right)  },$ we have the following result:

\begin{theorem}
\label{marta}Consider the equation
\begin{equation}
\partial_{t}u+(-1)^{m}\partial_{x}^{2m+1}u-|u|^{p-2}\partial_{x}u=0
\label{rita}%
\end{equation}
and assume that the Cauchy problem is globally well posed for (\ref{rita}) in
$H^{2m}.$Then (\ref{rita}) admits hylomorphic solitons $u_{\delta}$
$,\ \delta\in\left(  0,\delta_{\infty}\right)  $ provided that%
\begin{equation}
2<p<4m+2 \label{claretta}%
\end{equation}
Moreover, $u_{\delta}$ is a (weak) solution of the equation%
\begin{equation}
(-1)^{m}\partial_{x}^{2m+1}u-|u|^{p-2}\partial_{x}u=\lambda_{\delta}%
\partial_{x}u
\end{equation}
and
\begin{equation}
U_{\delta}(t,x)=u_{\delta}(x-\lambda_{\delta}t)
\end{equation}
solves (\ref{rita}) for suitable $\lambda_{\delta}$.
\end{theorem}

Observe that the Cauchy problem for (\ref{rita}) with $m=1$ is globally well
posed in $H^{2}$ (see\cite{We86}, \cite{kato}).

Theorem \ref{marta} generalizes a result of (\cite{befo2015}).

We consider also the generalized fractional nonlinear Schroedinger equation
(\ref{GNS}), since it is crucial in the study of (\ref{GBO}).

For $s=1$ and $W(u)=-\frac{1}{4}u^{4},$ (\ref{GNS}) reduces to the
Gross-Pitaevskii equation:
\begin{equation}
i\frac{\partial\psi}{\partial t}=-\frac{1}{2}\partial_{x}^{2}\psi-|\psi
|^{2}\psi,\ \ \psi(t,x)\in\mathbb{C}. \label{gpe}%
\end{equation}
Taking in (\ref{GNS}) $s=m$ and $W(r)=-\frac{\left\vert r\right\vert ^{p}}%
{p},$ we have the following theorem which generalizes some results relative to
(\ref{gpe}):

\begin{theorem}
\label{caterina}Consider the equation
\begin{equation}
i\frac{\partial\psi}{\partial t}=\frac{\left(  -1\right)  ^{m}}{2}\partial
_{x}^{2m}\psi-|\psi|^{p-2}\psi,\ \ \psi(t,x)\in\mathbb{C}, \label{nuovo}%
\end{equation}
and assume that the Cauchy problem is globally well posed in $H^{m}.$ Then
(\ref{nuovo}) admits hylomorphic solitons $u_{\delta}$ $,\ \delta\in\left(
0,\delta_{\infty}\right)  $ provided that (\ref{claretta}) holds. Moreover,
for suitable $\omega_{\delta},u_{\delta}$ is a (weak) solution of the equation%
\begin{equation}
\frac{\left(  -1\right)  ^{m}}{2}\partial_{x}^{2m}u-|u|^{p-2}u=\omega_{\delta
}u
\end{equation}
and
\begin{equation}
U_{\delta}(t,x)=u_{\delta}(x)e^{-i\omega_{\delta}t}%
\end{equation}
solves (\ref{nuovo}).
\end{theorem}

Observe that the Cauchy problem for (\ref{nuovo}) with $m=1$ is globally well
posed in $H^{1}$(see \cite{sulesulem} and its references)

\subsection{Notations\label{not}}

Let $\Omega$\ be a subset of $\mathbb{R}^{N}$: then

\begin{itemize}
\item $\mathcal{C}^{k}\left(  \mathbb{R}\right)  $ denotes the set of rael
functions which have continuous derivatives up to the order $k;$

\item $\mathcal{D}\left(  \mathbb{R}\right)  $ denotes the set of the
infinitely differentiable functions with compact support$;\ \mathcal{D}%
^{\prime}$ denotes the topological dual of $\mathcal{D}\left(  \mathbb{R}%
\right)  $, namely the set of distributions$;$

\item $H^{k}(a,b)$ is the closure in $L_{loc}^{1}$ of $\mathcal{D}\left(
\mathbb{R}\right)  $ with respect to the norm%
\[
\left\Vert u\right\Vert _{H^{k}(a,b)}^{2}=\int_{a}^{b}\left(  |D^{k}%
u(x)|^{2}+|u(x)|^{2}\right)  dx
\]

\item $\dot{H}^{s}$ is the closure in $L_{loc}^{1}$ of $\mathcal{D}\left(
\mathbb{R}\right)  $with respect to the norm%
\[
\left\Vert u\right\Vert _{\dot{H}^{s}}^{2}=\int|\xi|^{2s}|\hat{u}(\xi
)|^{2}d\xi
\]

\item $H^{s}\ $is the closure in $L_{loc}^{1}$ of $\mathcal{D}\left(
\mathbb{R}\right)  $ with respect to the norm%
\begin{align*}
\left\Vert u\right\Vert _{\dot{H}^{s}}^{2}  &  =\int\left(  1+|\xi
|^{2s}\right)  |\hat{u}(\xi)|^{2}d\xi\\
&  =\left\Vert u\right\Vert _{\dot{H}^{s}}^{2}+\left\Vert u\right\Vert
_{L^{2}}^{2}%
\end{align*}

\item $W^{s,2}(a,b),$ $s\in\left(  0,1\right)  ,$ is the closure in
$L_{loc}^{1}$ of $\mathcal{D}\left(  \mathbb{R}\right)  $ with respect to the
norm%
\[
\left\Vert u\right\Vert _{W^{s,2}(a,b)}^{2}=\int_{a}^{b}\int_{a}^{b}%
\frac{\left\vert u(x)-u(y)\right\vert ^{2}}{\left\vert x-y\right\vert ^{2s+1}%
}dxdy+\left\Vert u\right\Vert _{L^{2}(a,b)}^{2}%
\]

\item $W^{s,2}(a,b),s=k+\theta,$ $k\in\mathbb{N},\ \theta\in\left(
0,1\right)  ,$ is the closure in $L_{loc}^{1}$ of $\mathcal{D}\left(
\mathbb{R}\right)  $ with respect to the norm%
\[
\left\Vert u\right\Vert _{W^{k+\theta,2}(a,b)}^{2}=\left\Vert u\right\Vert
_{H^{k}(a,b)}^{2}+\left\Vert D^{k}u\right\Vert _{W^{\theta,2}(a,b)}^{2}%
\]

\end{itemize}

For the sake of the reader we will recall the properties of the Sobolev spaces
which we will use:%

\begin{equation}
\text{The norms of }\ \ W^{s,2}(\mathbb{R})\ \ \text{and\ \ }H^{s}%
\ \ \text{are equivalent} \label{mario}%
\end{equation}

For $a<b<c,$ we have that%
\begin{equation}
\left\Vert u\right\Vert _{W^{s,2}(a,b)}^{2}+\left\Vert u\right\Vert
_{W^{s,2}(,b,c)}^{2}\leq\left\Vert u\right\Vert _{W^{s,2}(a,c)}^{2}
\label{maria}%
\end{equation}

\section{Abstract theory\label{sws}}

In this section we construct an abstract functional framework which allows to
define solitary waves, solitons and hylomorphic solitons.

\subsection{Orbitally stable states and solitons\label{be}}

Let us consider the following field equation%
\begin{equation}
\frac{\partial\mathbf{u}}{\partial t}=\mathcal{A}(\mathbf{u}) \label{equa}%
\end{equation}
where $\mathbf{u}(t,\cdot)\in X,$ $X$ is a functional Hilbert space and
$\mathcal{A}:X\rightarrow Y$ is a differential operator.

We make the following assumptions:

\begin{itemize}
\item we have that%
\begin{equation}
X\subset L_{loc}^{1}\left(  \mathbb{R}^{N},V\right)  \label{lilla}%
\end{equation}
where $V$ is a vector space with norm $\left\vert \ \cdot\ \right\vert _{V}$
and which is called the internal parameters space

\item the dynamics has two constants of motion, the energy $E$ and the hylenic
charge $C;$ more precisely there are two continuous functionals on $X,$ which
are constant along the smooth solutions of (\ref{equa}); of course$\ $at this
level of abstraction the names energy and hylenic charge are conventional
\end{itemize}

For every $\mathbf{u}_{0}\in X$, the number $T(\mathbf{u}_{0})\geq0$ is
defined as the supremum of the $\tau\in\left[  0,+\infty\right)  $ such that
$\forall t\in\left[  0,\tau\right)  ,$ the solution $\mathbf{u}\left(
t,x\right)  $ of eq. \ref{equa} with initial value $\mathbf{u}_{0}\in X$
satisfies the following requests:

\begin{itemize}
\item $\mathbf{u}\left(  t,\cdot\right)  $ is the unique solution in $X;$

\item $\mathbf{u}\left(  \cdot,\mathbf{\cdot}\right)  $ is continuous in every
point $\left(  t,\mathbf{u}\right)  \in Z$ where%
\[
Z:=\left\{  (t,\mathbf{u})\in\mathbb{R}\times X\ |\ t\in\left[  0,T(\mathbf{u}%
)\right)  \right\}
\]

\item the functions $t\longmapsto E(\mathbf{u}\left(  t,\cdot\right)
)\ \ $and$\ \ t\longmapsto C(\mathbf{u}\left(  t,\cdot\right)  )$ are constant.
\end{itemize}

\bigskip

Also, we shall use the following notation: for every $\mathbf{u}_{0}\in X,$ we
set%
\begin{equation}
\gamma_{t}\mathbf{u}_{0}:=\mathbf{u}\left(  t,x\right)  ,\ t\in\left[
0,T(\mathbf{u}_{0})\right)  . \label{flusso}%
\end{equation}
where $\mathbf{u}\left(  t,x\right)  $ is the solution of eq. \ref{equa} with
initial value $\mathbf{u}_{0}.$ So that $\left(  X,\gamma\right)  $ defines a
dynamical system. Finally, we set%
\[
X_{0}=\left\{  \mathbf{u}\in X~|\ T(\mathbf{u})=+\infty\right\}
\]
In the case of KdV equation we have that $X=H^{1}$ (see \cite{befo2015}) and
$X_{0}=H^{2}\ $(see \cite{kato}, \cite{We86}).

\bigskip

Roughly speaking a \textit{soliton} is a localized state whose evolution
preserves this localization and which exhibits some form of stability so that
it has a particle-like behavior.\ To give a precise definition of
\textit{soliton} at this level of abstractness, we need to recall some well
known notions in the theory of dynamical systems.

We start with the definition of dynamical system suitable to our purposes:

\begin{definition}
A dynamical system is a couple $\left(  X,\gamma\right)  $ where $X$ is a
metric space and
\[
\gamma:Z\rightarrow X
\]
is a continuous map such that%
\begin{align*}
\gamma_{0}(\mathbf{u})  & =\mathbf{u}\\
\gamma_{t}\circ\gamma_{s}(\mathbf{u})  & =\gamma_{t+s}(\mathbf{u})
\end{align*}
provided that  $t+s\in\left[  0,T(\mathbf{u})\right)  .$
\end{definition}

\begin{definition}
A set $\Gamma\subset X$ is called \textit{invariant} if $\forall\mathbf{u}%
\in\Gamma,\forall t\in\left[  0,T(\mathbf{u})\right)  ,\ \gamma_{t}%
\mathbf{u}\in\Gamma.$
\end{definition}

\begin{definition}
An invariant set $\Gamma\subset X$ is called (orbitally) stable, if
$\forall\varepsilon>0,$ $\exists\delta>0,\;\forall\mathbf{u}\in X$,
\[
d(\mathbf{u},\Gamma)\leq\delta,
\]
implies that
\[
\forall t\in\left[  0,T(\mathbf{u})\right)  ,\text{ }d(\gamma_{t}%
\mathbf{u,}\Gamma)\leq\varepsilon.
\]

\end{definition}

\begin{remark}
Notice that, under reasonable assumptions, if $\Gamma$ is stable, then%
\[
\exists\delta>0,\;d(\mathbf{u},\Gamma)\leq\delta\Rightarrow T(\mathbf{u}%
)=+\infty
\]

\end{remark}

Now we can give our definition of soliton:

\begin{definition}
\label{dos} A soliton is a state $\mathbf{u}\in\Gamma\subset X_{0},$ where
$\Gamma$ satisfies the following properties:

\begin{itemize}
\item (i) $\Gamma$ is an invariant, stable set,

\item (ii) $\Gamma$ is "compact up translations", namely for any sequence
$\mathbf{u}_{n}(x)\in\Gamma\ $there is a subsequence $\mathbf{u}_{n_{k}}$ and
a sequence $\tau_{k}\in\mathbb{R}^{n}$ such that $\mathbf{u}_{n_{k}}%
(x-\tau_{k})$ is convergent.
\end{itemize}
\end{definition}

\begin{remark}
The above definition needs some explanation. For simplicity, we assume that
$\Gamma$ is a manifold (actually, it is possible to prove that this is the
generic case if the problem is formulated in a suitable function space). Then
(ii) implies that $\Gamma$ is finite dimensional. Since $\Gamma$ is invariant,
$\mathbf{u}_{0}\in\Gamma\Rightarrow\gamma_{t}\mathbf{u}_{0}\in\Gamma$ for
every time. Thus, since $\Gamma$ is finite dimensional, the evolution of
$\mathbf{u}_{0}$ is described by a finite number of parameters. By the
stability of $\Gamma$, a small perturbation of $\mathbf{u}_{0}$ remains close
to $\Gamma.$ However, in this case, its evolution depends on an infinite
number of parameters. Thus, this system appears as a finite dimensional system
with a small perturbation.
\end{remark}

\subsection{An abstract theorem}

We recall that our dynamical system $\left(  X,\gamma\right)  $ has two
constants of motion: the energy $E$ and the hylenic charge $C.$

\begin{definition}
\label{tdc} A soliton $\mathbf{u}_{0}\in X$ is called \textbf{hylomorphic }if
the set $\Gamma$ (given by Def. \ref{dos}) has the following structure%
\begin{equation}
\Gamma=\Gamma\left(  e_{0},c_{0}\right)  =\left\{  \mathbf{u}\in
X\ |\ E(\mathbf{u})=e_{0},\ \left\vert C(\mathbf{u})\right\vert =c_{0}%
\right\}  \label{plis}%
\end{equation}
where%
\begin{equation}
e_{0}=\min\left\{  E(\mathbf{u})\ |\ \left\vert C(\mathbf{u})\right\vert
=c_{0}\right\}  . \label{minbis}%
\end{equation}

\end{definition}

Notice that, by (\ref{minbis}), we have that a hylomorphic soliton
$\mathbf{u}_{0}$ minimizes the energy on%
\begin{equation}
\mathfrak{M}_{c_{0}}=\left\{  \mathbf{u}\in X\ |\ \left\vert C(\mathbf{u}%
)\right\vert =c_{0}\right\}  .
\end{equation}
If $\mathfrak{M}_{c_{0}}$ is a manifold and $E$ and $C$ are differentiable,
then $\mathbf{u}_{0}$ satisfies the following nonlinear eigenvalue problem:%
\[
E^{\prime}(\mathbf{u}_{0})=\lambda C^{\prime}(\mathbf{u}_{0}).
\]

We assume that $E$ and $C$ satisfy the following assumptions:

\begin{itemize}
\item \textit{(EC-1) \textbf{(Value at 0)} }$E,$ $C$ \textit{are} $C^{1}$,
bounded on bounded sets functionals and such that\textit{\ }%
\[
E(0)=0,\ C(0)=0;\ E^{\prime}(0)=0;\ C^{\prime}(0)=0.
\]

\item \textit{(EC-2) \textbf{(Invariance)\ }}$E$ \textit{and }$C$\textit{\ are
}$\tau-$\textit{invariant i.e. invariant under space translations..}

\item \textit{(EC-3)\textbf{(Coercivity)} we assume that }$C(\mathbf{u}%
)>0\Leftrightarrow\mathbf{u}\neq0$ and\textit{ that} \textit{there exists
}$a\geq0$ \textit{and }$\beta>1$\textit{\ such that}

\begin{itemize}
\item (i) $E(\mathbf{u})+aC(\mathbf{u})^{\beta}\geq0;$

\item (ii) \textit{if }$\left\Vert \mathbf{u}\right\Vert \rightarrow\infty
,\ $\textit{then} $E(\mathbf{u})+aC(\mathbf{u})^{\beta}\rightarrow\infty;$

\item (iii) \textit{for any} \textit{bounded sequence }$\mathbf{u}_{n}$
\textit{in }$X$ \textit{such that} $E(\mathbf{u}_{n})+aC(\mathbf{u}%
_{n})^{\beta}\rightarrow0,\ $\textit{we have that }$\mathbf{u}_{n}%
\rightarrow0.$
\end{itemize}

\item \textit{(EC-4)\textbf{(Splitting property)} }$E$\textit{\ and }%
$C$\textit{\ satisfy the splitting property.~}We say that a functional $F$ on
$X$ has the splitting property if given a sequence $\mathbf{u}_{n}%
=\mathbf{u}+\mathbf{w}_{n}$ in $X$ such that $\mathbf{w}_{n}$ converges weakly
to $0$, we have that
\begin{equation}
F(\mathbf{u}_{n})=F(\mathbf{u})+F(\mathbf{w}_{n})+o(1).
\end{equation}

\end{itemize}

We need also the following definition:

\begin{definition}
\label{sega}A sequence $\psi_{n}\ $is called vanishing if

\begin{itemize}
\item \textit{for }any sequence $\left\{  x_{n}\right\}  \subset\mathbb{R}%
^{N}$ the translated sequence $\left\{  \psi_{n}\mathbf{(\cdot-}%
x_{n}\mathbf{)}\right\}  $ converges weakly in $X$ to $0.$
\end{itemize}
\end{definition}

\textbf{Example.} If $\psi_{n}$ converges strongly to zero then it is a
vanishing sequence, but the converse is not true; for example the sequence%
\[
\psi_{n}=e^{-x^{2}}\sin nx
\]
in $L^{2}$; it does not converge strongly to 0 and it is vanishing.

\bigskip

We set
\begin{equation}
\Lambda\left(  \mathbf{u}\right)  :=\frac{E\left(  \mathbf{u}\right)
}{\left\vert C\left(  \mathbf{u}\right)  \right\vert }. \label{lambda}%
\end{equation}

Since $E$ and $C$ are constants of motion, also $\Lambda$ is a constant of
motion; it will be called \textbf{hylenic ratio} and, as we will see it will
play a central role in this theory. Finally we set
\begin{equation}
\Lambda_{0}:=\ \inf\left\{  \lim\inf\ \Lambda(\mathbf{u}_{n})\ |\ \mathbf{u}%
_{n}\ \text{is a vanishing sequence}\right\}  \label{hylo}%
\end{equation}

\begin{theorem}
\label{astra1}Assume that $E\ $and $C$ satisfy (EC-1),(EC-2),\textit{(EC-3).
Moreover assume that }the following condition
\begin{equation}
\underset{u}{\inf}\Lambda(\mathbf{u})<\Lambda_{0}. \label{hh}%
\end{equation}
is satisfied. Then for every $\delta\in\left(  0,\delta_{\infty}\right)  ,$
$\delta_{\infty}>0,$ there exist $c_{\delta}>0$ and $\Gamma_{\delta}$ $\subset
X_{0}$ satisfying i),ii) of definition \ref{dos} and such that any soliton
$\mathbf{u}_{\delta}\in\Gamma_{\delta}$ is hylomorphic, i.e. it minimizes the
energy on the manifold%
\[
\mathfrak{M}_{c_{\delta}}=\left\{  \mathbf{u}\in X\ |\ C(\mathbf{u}%
)=c_{\delta}\right\}  .
\]
Moreover if $\delta_{1}<\delta_{2}$ we have that $c_{\delta_{1}}>c_{\delta
_{2}})$.\ 
\end{theorem}

The proof of this theorem is in \cite{befolib} Th. 34 pag. 39.

\bigskip

The inequality (\ref{hh}) plays a crucial role in this theory; we will refer
to it as to the \textbf{hylomorphy condition}.

\section{The nonlinear fractional Schr\"{o}dinger equation\label{NSE}}

\bigskip

\subsection{Main results}

The solitons for eq. (\ref{GBO}), as we will see, are related to the solitons
of the Fractional Nonlinear Schr\"{o}dinger equation (\ref{GNS}).

Here we shall use a method to prove the existence of hylomorphic solitons for
(\ref{GBO}) similar to the one presented in \cite{befolib} (see also
\cite{befolak} and \cite{befobrez}). In this section we will resume this method.

The Fractional Nonlinear Schr\"{o}dinger equation is given by
\begin{equation}
i\frac{\partial\psi}{\partial t}=\frac{1}{2}D^{2s}\psi+\frac{1}{2}W^{\prime
}(\psi);\ \ s>0, \label{NSV}%
\end{equation}
where $\psi:\mathbb{R\times R}\rightarrow\mathbb{C}$, $D$ is defined by
(\ref{df}), $W:\mathbb{C\rightarrow R}$ and
\begin{equation}
W^{\prime}(\psi)=\frac{\partial W}{\partial\psi_{1}}+i\frac{\partial
W}{\partial\psi_{2}}. \label{w'}%
\end{equation}

We assume that $W$ depends only on $\left\vert \psi\right\vert $, namely
\[
W(\psi)=F(\left\vert \psi\right\vert )\ \text{and so\ }W^{\prime}%
(\psi)=F^{\prime}(\left\vert \psi\right\vert )\frac{\psi}{\left\vert
\psi\right\vert }.
\]
for some smooth function $F:\left[  0,\infty\right)  \rightarrow\mathbb{R}.$

\begin{proposition}
Let us consider the dynamical system $(H^{s},\gamma)$ related to eq.
(\ref{NSV}). Then the energy%
\begin{equation}
E=\int\left(  \frac{1}{2}\left\vert D^{s}\psi\right\vert ^{2}+W(\psi)\right)
dx\label{ghj}%
\end{equation}
and the charge%
\begin{equation}
C=\int\left\vert \psi\right\vert ^{2}dx\label{chery}%
\end{equation}
are constant along a smooth solution $\psi$ which decay in space sufficiently fast.
\end{proposition}

\textbf{Proof:} Let $\psi\ $be a smooth solution of (\ref{NSV}) with initial
condition $\psi_{0}$. Then, we have
\begin{align*}
\frac{d}{dt}E(\psi(t))  &  =\operatorname{Re}\int\left(  D^{s}\psi
\ \partial_{t}D^{s}\overline{\psi}+W^{\prime}(\psi)\partial_{t}\overline{\psi
}\right)  dx\\
&  =\operatorname{Re}\int\left(  D^{2s}\psi+W^{\prime}(\psi)\right)
\overline{\partial_{t}\psi}dx\\
&  =\operatorname{Re}\int\left(  D^{2s}\psi+W^{\prime}(\psi)\right)
\overline{\left(  \frac{1}{2}iD^{2s}\psi+\frac{1}{2}iW^{\prime}(\psi)\right)
}\\
&  =-\operatorname{Re}\left[  \frac{1}{2}i\int\left(  D^{2s}\psi+W^{\prime
}(\psi)\right)  \overline{\left(  D^{2s}\psi+W^{\prime}(\psi)\right)  }\right]
\\
&  =-\operatorname{Re}\left[  \frac{1}{2}i\int\left\vert D^{2s}\psi+W^{\prime
}(\psi)\right\vert ^{2}dx\right]  =0.
\end{align*}

Moreover%
\begin{align*}
\frac{d}{dt}\int\left\vert \psi\right\vert ^{2}dx  &  =2\operatorname{Re}%
\int\left(  \psi\ \partial_{t}\overline{\psi}\right)  dx\\
&  =-2\operatorname{Re}i\int\left(  \psi\ \overline{\left(  D^{2s}%
\psi+W^{\prime}(\psi)\right)  }\right)  dx=0
\end{align*}

$\square$

\bigskip

We make the following assumptions on the function $W:$
\begin{equation}
W(0)=W^{\prime}(0)=0. \tag{W-0}\label{Wa}%
\end{equation}

Set
\begin{equation}
W(r)=E_{0}r^{2}+N(r),\ \ E_{0}=\frac{1}{2}W^{\prime\prime}(0) \label{W}%
\end{equation}
and assume that
\begin{equation}
\exists r_{0}\in\mathbb{R}^{+}\ \text{such that\ }N(r_{0})<0. \tag{W-1}%
\label{W1}%
\end{equation}
There exist $q_{1}\leq q_{2}$ in $(2,+\infty),$ s. t.%
\begin{equation}
|N^{\prime}(r)|\leq c_{1}\left(  r^{q_{1}-1}+r^{q_{2}-1}\right)
\tag{W-2}\label{Wp}%
\end{equation}
Moreover assume that, $\exists p\in\left(  2,4s+2\right)  ,$with $s\geq
\frac{1}{2},$s.t.$\ $
\begin{equation}
N(r)\geq-cr^{p},\text{ }c\geq0,\ \text{ for }r\text{ large} \tag{W-3}%
\label{W0}%
\end{equation}

We can apply the abstract theory of section \ref{sws} taking $X=H^{s}$:

\begin{theorem}
\label{solitoni} Let $W$ satisfy (\ref{Wa}),...,(\ref{W0}), and
\begin{equation}
\frac{1}{2}W^{\prime\prime}(0)=E_{0}>0\label{Wb}%
\end{equation}
Then equation (\ref{NSV}) admits a family of hylomorphic solitons $u_{\delta
},\ \delta\in\left(  0,\delta_{\infty}\right)  ,\ u_{\delta}\in H^{2s}.$
\end{theorem}

\textbf{Proof: }The proof of this theorem will be given in the next section.

$\square$

The following theorem gives more information on the structure of the solitons
and of their dynamics; moreover assumption (\ref{Wb}) is avoided.

\begin{theorem}
\label{marina}Assume that all the hypotheses of Theorem \ref{solitoni} hold
with exception of assumption (\ref{Wb}). Then equation (\ref{NSV}) admits
hylomorphic solitons $u_{\delta}\in H^{2s},\ \delta\in\left(  0,\delta
_{\infty}\right)  .$ Moreover, $u_{\delta}$ is a (weak) solution of the
equation%
\begin{equation}
\frac{1}{2}D^{2s}u+\frac{1}{2}W^{\prime}(u)=\omega u\label{brb}%
\end{equation}
and
\begin{equation}
\psi_{\delta}\left(  t,x\right)  :=u_{\delta}(x)e^{-i\omega t}\label{lula2}%
\end{equation}
solves (\ref{NSV}).
\end{theorem}

\textbf{Proof}. First let us assume (\ref{Wb}). By Theorem \ref{solitoni}
(\ref{NSV}), admits a family of hylomorphic solitons $u_{\delta},\ \delta
\in\left(  0,\delta_{\infty}\right)  ,\ u_{\delta}\in H^{2s}.$.

Let $u_{\delta}$ be a hylomorphic soliton, then it is a minimizer of the
energy $E$ defined in (\ref{ghj}) on $\mathfrak{M}_{c}.$ Then we get
\begin{equation}
E^{\prime}(u_{\delta})=\omega C^{\prime}(u_{\delta}) \label{bekime}%
\end{equation}
where $\omega$ is a Lagrange multiplier. Clearly (\ref{bekime}) gives
\[
\frac{1}{2}D^{2s}u_{\delta}+\frac{1}{2}W^{\prime}(u_{\delta})=\omega
u_{\delta}%
\]
which implies that $\psi=u_{\delta}e^{-i\omega t}$ solves (\ref{NSV}). Observe
that, since $u_{\delta}\in H^{s}$ solves (\ref{bekime}), by elliptic
reguralization we have $u_{\delta}\in H^{2s}.$ It remains to show that the
same result holds even when (\ref{Wb}) is violated. So assume that
\[
\frac{1}{2}W^{\prime\prime}(0)=E_{0}\leq0;
\]
we can reduce the problem to the case (\ref{Wb}). To do this, we replace
$W(r)$ with
\[
W_{1}(r)=W(r)+\frac{1}{2}\left(  1-E_{0}\right)  r^{2}%
\]
So
\[
\frac{1}{2}W_{1}^{\prime\prime}(0)=1>0.
\]
and (\ref{Wb}) is satisfied by $W_{1}$. Then we can apply the previous
considerations. So there exists a hylomorphic soliton $u_{1}$ and $\omega_{1}$
s.t.%
\[
\psi_{1}=u_{1}e^{-i\omega_{1}t}.
\]
solves the equation
\begin{equation}
i\frac{\partial\psi_{1}}{\partial t}=\frac{1}{2}D^{2s}\psi_{1}+\frac{1}%
{2}W_{1}^{\prime}(\psi_{1}), \label{prooo}%
\end{equation}
It can be easily seen that $\psi=\psi_{1}\left(  t,x\right)  e^{i\frac{\left(
E_{0}-1\right)  }{2}t}$ is a solution of (\ref{NSV})

$\square$

\bigskip

\begin{remark}
Th. \ref{caterina} is a particular case of the above theorem when $s$ in
(\ref{NSV}) is a positive integer.
\end{remark}

\subsection{Proof of Th. \ref{solitoni}}

In this section, we will prove Theorem \ref{solitoni} by using Th.\ref{astra1}
with $X=H^{s},\ \ E$ and $C$ as in (\ref{ghj}) and (\ref{chery}) and
\begin{equation}
\Lambda\left(  \psi\right)  :=\frac{E\left(  \psi\right)  }{\left\vert
C\left(  \psi\right)  \right\vert }=\frac{\int\left(  \frac{1}{2}\left\vert
D^{s}\psi\right\vert ^{2}+W(\psi)\right)  dx}{\int\left\vert \psi\right\vert
^{2}dx}.
\end{equation}
Let us first prove the Coercivity assumption (\textit{EC-3)};

\begin{lemma}
\label{EC3}If $W$ satisfies (\ref{Wa}),...(\ref{W0}) then assumpion (EC-3) holds
\end{lemma}

\textbf{Proof:} by (\ref{W}) and (\ref{W0}), we have
\begin{align}
E(\mathbf{u})+aC(\mathbf{u})^{\beta}  &  =\int\left(  \frac{1}{2}\left\vert
D^{s}\psi\right\vert ^{2}+W(\psi)\right)  dx+a\left(  \int\left\vert
\psi\right\vert ^{2}dx\right)  ^{\beta}\nonumber\\
&  \geq\frac{1}{2}\left\Vert u\right\Vert _{\dot{H}^{s}}^{2}+E_{0}\left\Vert
u\right\Vert ^{2}-b\left\Vert u\right\Vert _{L^{p}}^{p}+a\left\Vert
u\right\Vert _{L^{2}}^{2\beta}.\nonumber\\
&  \geq\frac{1}{2}\left\Vert u\right\Vert _{\dot{H}^{s}}^{2}-b\left\Vert
u\right\Vert _{L^{p}}^{p}+a\left\Vert u\right\Vert _{L^{2}}^{2\beta}.
\label{anna}%
\end{align}
where $b$ is a suitable constant and $a,\beta$ will be choosen later.

The Gagliardo Nirenberg inequalities in our case take the following form:%
\[
\left\Vert u\right\Vert _{L^{p}}\leq c\left\Vert u\right\Vert _{\dot{H}^{s}%
}^{\theta}\left\Vert u\right\Vert _{L^{2}}^{1-\theta}%
\]
provided that
\[
\frac{1}{p}=\frac{1-\theta}{2}+\theta\left(  \frac{1}{2}-s\right)
\ \text{and}\ \theta\in\left(  0,1\right)
\]
namely$\ $%
\begin{equation}
\theta=\frac{1}{s}\left(  \frac{1}{2}-\frac{1}{p}\right)  \ \text{and}%
\ \theta\in\left(  0,1\right)  \label{paula}%
\end{equation}
Notice that for $s\geq\frac{1}{2}$ and $p>2,$ the above conditions are satisfied.

Then,
\begin{align*}
\left\Vert u\right\Vert _{L^{p}}^{p}  &  \leq c^{p}\left\Vert u\right\Vert
_{\dot{H}^{s}}^{p\theta}\left\Vert u\right\Vert _{L^{2}}^{p-p\theta}\\
&  =\frac{1}{2bp\theta}\left(  \left\Vert u\right\Vert _{\dot{H}^{s}}%
^{2}\right)  ^{\left(  p\theta\right)  /2}\left(  c_{1}\left\Vert u\right\Vert
_{L^{2}}^{p-p\theta}\right)
\end{align*}

By (\ref{W}) $p<4s+2,$ and hence $\frac{2}{p\theta}=\frac{4s}{p-2}>1.$ We now
use Young's inequality:%
\begin{align}
\left\Vert u\right\Vert _{L^{p}}^{p}  &  \leq\frac{1}{2bp\theta}\left(
\frac{p\theta}{2}\left[  \left(  \left\Vert u\right\Vert _{\dot{H}^{s}}%
^{2}\right)  ^{\left(  p\theta\right)  /2}\right]  ^{2/\left(  p\theta\right)
}+c_{2}\left(  \left\Vert u\right\Vert _{L^{2}}^{p-p\theta}\right)  ^{\left(
2/\left(  p\theta\right)  \right)  ^{\prime}}\right) \label{lisa}\\
&  =\frac{1}{4b}\left\Vert u\right\Vert _{\dot{H}^{s}}^{2}+\frac{c_{2}%
}{2bp\theta}\left\Vert u\right\Vert _{L^{2}}^{2\beta}\nonumber
\end{align}

where $\beta=\left(  \frac{p-p\theta}{2}\right)  \left(  \frac{2}{p\theta
}\right)  ^{\prime}.$ We want to show that $\beta>1.$ We have that%
\[
\beta=\left(  1-\theta\right)  \frac{p}{2-p\theta}%
\]
and using (\ref{paula}), since $\beta>0\ $and $\ p>2,$ we get%
\[
\beta=\frac{\left(  2ps+2\right)  -p}{(4s+2)-p}>1
\]

So, by the above inequalities (\ref{lisa}) and (\ref{anna}),
\begin{align}
E(\mathbf{u})+aC(\mathbf{u})^{\beta}  &  \geq\frac{1}{2}\left\Vert
u\right\Vert _{\dot{H}^{s}}^{2}-b\left[  \frac{1}{4b}\left\Vert u\right\Vert
_{\dot{H}^{s}}^{2}+\frac{c_{2}}{b}\left\Vert u\right\Vert _{L^{2}}^{2\beta
}\right]  +a\left\Vert u\right\Vert _{L^{2}}^{2\beta}.\\
&  =\frac{1}{4}\left\Vert u\right\Vert _{\dot{H}^{s}}^{2}+\left(
a-c_{3}\right)  \left\Vert u\right\Vert _{L^{2}}^{2\beta}%
\end{align}

If $a\geq c_{3},$ then all the assumptions \textit{(EC-3) }are easily verified.

$\square$

\bigskip

\begin{lemma}
\label{EC4}If $W$ satisfies (\ref{Wa}),...(\ref{W0}) then the splitting
property (EC-4) holds.
\end{lemma}

\textbf{Proof:} See the proof Lemma 5.3 at pg. 74 in \cite{befolib}.

$\square$

\bigskip

Next we will verify that the hylomorphy condition (\ref{hh}) is satisfied. The
following lemma, which is in the same spirit of some compactness results in
\cite{lieb}, \cite{bece} and \cite{brezis}, plays a fondamental role in
proving (\ref{hh}):

\begin{lemma}
\label{nonvanishing2}\textit{\textbf{\ }} Let $\psi_{n}\ $be a vanishing
sequence in $H^{s}$, $s\geq\frac{1}{2}$ (see Def. \ref{sega}); then for any
$p>2$ we have $\left\Vert \psi_{n}\right\Vert _{L^{p}}\rightarrow0.$
\end{lemma}

\textbf{Proof. }Now let $\left\{  \psi_{n}\right\}  \subset H^{s}$ be a
vanishing sequence and prove that $\left\Vert \psi_{n}\right\Vert _{L^{p}%
}\rightarrow0$. Arguing by contradiction, assume that, up to a subsequence,
$\left\Vert \psi_{n}\right\Vert _{L^{p}}\geq a>0.$ Since $\psi_{n}$ is
vanishing, there exists $M>0$ such that $\left\Vert \psi_{n}\right\Vert
_{H^{1}}^{2}\leq M.$ Then, if $L$ is the constant for the Sobolev embedding
$H^{s}\subset L^{p}\left(  j,j+1\right)  ,$ we have
\begin{align*}
0  &  <a^{p}\leq\int\left\vert \psi_{n}\right\vert ^{p}=\sum_{j}\int_{j}%
^{j+1}\left\vert \psi_{n}\right\vert ^{p}=\sum_{j}\left\Vert \psi
_{n}\right\Vert _{L^{p}\left(  j,j+1\right)  }^{p-2}\left\Vert \psi
_{n}\right\Vert _{L^{p}\left(  j,j+1\right)  }^{2}\\
&  \leq\ \left(  \underset{j}{\sup}\left\Vert \psi_{n}\right\Vert
_{L^{p}\left(  j,j+1\right)  }^{p-2}\right)  \cdot\sum_{j}\left\Vert \psi
_{n}\right\Vert _{L^{p}\left(  j,j+1\right)  }^{2}\\
&  \leq\ L\left(  \underset{j}{\sup}\left\Vert \psi_{n}\right\Vert
_{L^{p}\left(  j,j+1\right)  }^{p-2}\right)  \cdot\sum_{j}\left\Vert \psi
_{n}\right\Vert _{W^{2,s}\left(  j,j+1\right)  }^{2}\ \text{(by (\ref{maria}%
))}\\
&  \leq L\left(  \underset{j}{\sup}\left\Vert \psi_{n}\right\Vert
_{L^{p}\left(  \left[  j,j+1\right]  \right)  }^{p-2}\right)  \left\Vert
\psi_{n}\right\Vert _{W^{2,s}(\mathbb{R)}}^{2}\ \ \text{(by (\ref{mario}))}\\
&  \leq LM_{1}\left(  \underset{j}{\sup}\left\Vert \psi_{n}\right\Vert
_{L^{p}\left(  \left[  j,j+1\right]  \right)  }^{p-2}\right)  \left\Vert
\psi_{n}\right\Vert _{H^{s}}^{2}\leq LM_{1}M\left(  \underset{j}{\sup
}\left\Vert \psi_{n}\right\Vert _{L^{p}\left(  \left[  j,j+1\right]  \right)
}^{p-2}\right)  .
\end{align*}
Then%
\[
\left(  \underset{j}{\sup}\left\Vert \psi_{n}\right\Vert _{L^{p}\left(
\left[  j,j+1\right]  \right)  }\right)  \geq\left(  \frac{a^{p}}{LM}\right)
^{1/(p-2)}%
\]

Then, for any $n,$ there exists $j_{n}\in\mathbb{Z}$ such that
\begin{equation}
\left\Vert \psi_{n}\right\Vert _{L^{p}\left(  j_{n},j_{n}+1\right)  }%
\geq\alpha>0. \label{caca}%
\end{equation}
Then, we easily have
\begin{equation}
\left\Vert \psi_{n}\left(  \cdot-j_{n}\right)  \right\Vert _{L^{p}%
(0,1)}=\left\Vert \psi_{n}\right\Vert _{L^{p}\left(  j_{n},j_{n}+1\right)
}\geq\alpha>0. \label{chicco}%
\end{equation}

Since $\psi_{n}$ is bounded, also $\psi_{n}\left(  \cdot-j_{n}\right)  $ is
bounded (in $H^{s}).$ Then we have, up to a subsequence, that $\psi_{n}\left(
\cdot-j_{n}\right)  \rightharpoonup\psi_{0}$ weakly in $H^{s}$ and hence
strongly in $L^{p}(0,1)$. By (\ref{chicco}), $\psi_{0}\neq0$ and this
contradicts the fact that $\psi_{n}$ is vanishing.

$\square$

\begin{lemma}
\label{preparatorio}If the assumptions of Theorem \ref{solitoni} are
satisfied, we have%
\[
\underset{\psi\in H^{s},\left\Vert \psi\right\Vert _{L^{p}}\rightarrow0}%
{\lim\inf}\Lambda(\psi)\geq E_{0}%
\]

\end{lemma}

\textbf{Proof. }Clearly%

\begin{align*}
\underset{\psi\in H^{s},\left\Vert \psi\right\Vert _{L^{p}}\rightarrow0}%
{\lim\inf}\Lambda(\psi)  &  =\ \underset{\psi\in H^{s},\left\Vert
\psi\right\Vert _{L^{p}}=1,\varepsilon\rightarrow0}{\lim\inf}\frac
{E(\varepsilon\psi)}{C(\varepsilon\psi)}\\
&  =\underset{\psi\in H^{s},\left\Vert \psi\right\Vert _{L^{p}}=1,}{\inf
}\left(  \frac{\int\left(  \frac{1}{2}\left\vert D^{s}\psi\right\vert
^{2}+E_{0}\left\vert \psi\right\vert ^{2}\right)  dx}{\int\left\vert
\psi\right\vert ^{2}}\right)  +\underset{\psi\in H^{s},\left\Vert
\psi\right\Vert _{L^{p}}=1,\varepsilon\rightarrow0}{\lim\inf}\frac{\int
N(\varepsilon\psi)}{\varepsilon^{2}\int\left\vert \psi\right\vert ^{2}}\\
&  \geq E_{0}+\underset{\psi\in H^{s},\left\Vert \psi\right\Vert _{L^{p}%
}=1,\varepsilon\rightarrow0}{\lim\inf}\frac{\int N(\varepsilon\psi
)}{\varepsilon^{2}\int\left\vert \psi\right\vert ^{2}}%
\end{align*}
So the proof of Lemma will be achieved if we show that%
\begin{equation}
\underset{\psi\in H^{s},\left\Vert \psi\right\Vert _{L^{p}}=1,\varepsilon
\rightarrow0}{\lim\inf}\frac{\int N(\varepsilon\psi)}{\varepsilon^{2}%
\int\left\vert \psi\right\vert ^{2}}=0. \label{resto}%
\end{equation}
By assumptions (\ref{Wp}) and (\ref{W0}) we have%
\begin{equation}
-cr^{p}\leq N(r)\leq c_{1}\left(  r^{q_{1}-1}+r^{q_{2}-1}\right)  .
\label{zerobis}%
\end{equation}

Then by (\ref{zerobis}) we have%
\[
-\underset{\left\Vert \psi\right\Vert _{L^{p}}=1}{\inf}\frac{c\int\left\vert
\varepsilon\psi\right\vert ^{p}}{\varepsilon^{2}\int\left\vert \psi\right\vert
^{2}}\leq\underset{\left\Vert \psi\right\Vert _{L^{p}}=1}{\inf}\frac{\int
N(\varepsilon\psi)}{\varepsilon^{2}\int\left\vert \psi\right\vert ^{2}}%
\leq\underset{\left\Vert \psi\right\Vert _{L^{p}}=1}{\inf}\frac{\int
c_{1}\left(  \left\vert \varepsilon\psi\right\vert ^{q_{1}-1}+\left\vert
\varepsilon\psi\right\vert ^{q_{2}-1}\right)  }{\varepsilon^{2}\int\left\vert
\psi\right\vert ^{2}}%
\]

\begin{equation}
-cA\varepsilon^{p-2}\leq\underset{\left\Vert \psi\right\Vert _{L^{p}}=1}{\inf
}\frac{\int N(\varepsilon\psi)}{\varepsilon^{2}\int\left\vert \psi\right\vert
^{2}}\leq c_{1}B(\varepsilon^{q_{1}-1}+\varepsilon^{q_{2}-1}) \label{uno}%
\end{equation}
where%
\[
A=\underset{\psi\in H^{s}\text{ }\left\Vert \psi\right\Vert _{L^{p}}=1}{\inf
}\frac{\int\left\vert \psi\right\vert ^{p}}{\int\left\vert \psi\right\vert
^{2}},\text{ }B=\underset{\psi\in H^{s}\text{ }\left\Vert \psi\right\Vert
_{L^{p}}=1}{\inf}\frac{\int\left(  \left\vert \psi\right\vert ^{q_{1}%
-1}+\left\vert \psi\right\vert ^{q_{2}-1}\right)  }{\int\left\vert
\psi\right\vert ^{2}}.
\]
By (\ref{uno}) we easily get (\ref{resto}).

$\square$\bigskip

\begin{corollary}
\label{interm} If the assumptions of Theorem \ref{solitoni} are satisfied,
then%
\[
E_{0}\leq\Lambda_{0}%
\]

\end{corollary}

\textbf{Proof.} By Lemma \ref{nonvanishing2} and Lemma \ref{preparatorio}%

\begin{align*}
\Lambda_{0}\  &  =\inf\left\{  \lim\inf\ \Lambda(\mathbf{u}_{n}%
)\ |\ \mathbf{u}_{n}\ \text{is a vanishing sequence}\right\} \\
&  \geq\ \underset{\left\Vert \psi\right\Vert _{L^{p}}\rightarrow0}{\lim\inf
}\ \Lambda(\psi)\geq E_{0}%
\end{align*}

$\square$

Finally we can prove that the hylomorphy condition is satisfied.

\begin{lemma}
\label{seguente2}If the assumptions of Theorem \ref{solitoni} are satisfied,
then the hylomorphy condition (\ref{hh}) holds, namely we have
\[
\underset{\psi\in H^{s}}{\inf}\Lambda(\psi)<\Lambda_{0}%
\]

\end{lemma}

\textbf{Proof. } We need to construct a function $u\in H^{s}$ such that
$\Lambda(u)<\Lambda_{0}.$ Such a function can be constructed as follows. Let
$u_{R}\geq0$ be a $C^{\infty}$ function such that
\[
u_{R}=\left\{
\begin{array}
[c]{cc}%
s_{0} & if\;\;|x|<R\\
0 & if\;\;|x|>R+1
\end{array}
.\right.
\]
and there exists a constant $C$ such that, for any integer $m\geqslant s\ $we
have
\[
\left\vert D^{m}u_{R}(x)\right\vert \leq C\ \
\]
There is a constant $c$ such that%
\[
\left\Vert u\right\Vert _{H^{m}}^{2}\leq c\left[  \int\left\vert
D^{m}u(x)\right\vert ^{2}dx+%
{\displaystyle\int}
\left\vert u(x)\right\vert ^{2}dx\right]
\]
then
\begin{align*}
\int\left\vert D^{s}u_{R}(x)\right\vert ^{2}dx  &  \leq c_{1}\left\Vert
u\right\Vert _{H^{m}}^{2}\leq cc_{1}\left[  \int\left\vert D^{m}%
u_{R}(x)\right\vert ^{2}dx+%
{\displaystyle\int}
\left\vert u_{R}(x)\right\vert ^{2}dx\right] \\
&  \leq2C^{2}cc_{1}+2s_{0}^{2}(R+1)\leq c_{2}%
\end{align*}
Moreover,
\[
2s_{0}^{2}R\leq\int\left\vert u_{R}\right\vert ^{2}dx\leq2s_{0}^{2}\left(
R+1\right)
\]
then%

\begin{equation}
\frac{\int\left[  \frac{1}{2}\left\vert D^{s}u_{R}\right\vert ^{2}+E_{0}%
u_{R}^{2}\right]  dx}{\int u_{R}^{2}}\leq E_{0}+O\left(  \frac{1}{R}\right)  .
\label{pinco}%
\end{equation}
Moreover
\[
\int N(u_{R})dx=2RN(s_{0})+\int_{R}^{R+1}N(u_{R})dx+\int_{-R-1}^{-R}%
N(u_{R})dx.
\]
So%
\begin{align}
\frac{\int N(u_{R})dx}{\int u_{R}^{2}}  &  \leq\frac{2RN(s_{0})+c_{3}}{\int
u_{R}^{2}}\leq(\text{ since }N(s_{0})<0)\label{palle}\\
&  \leq\frac{2RN(s_{0})}{2s_{0}^{2}(R+1)}+\frac{c_{3}}{2s_{0}^{2}R}%
=\frac{N(s_{0})}{s_{0}^{2}}+O\left(  \frac{1}{R}\right)  .\nonumber
\end{align}

Then, by (\ref{pinco}) e (\ref{palle}) we get%
\begin{align}
\Lambda(u_{R})  &  =\frac{\int\left(  \frac{1}{2}\left\vert D^{s}%
u_{R}\right\vert ^{2}+W(u_{R})\right)  dx}{\int u_{R}^{2}dx}\\
&  =\ \frac{\int\left(  \frac{1}{2}\left\vert D^{s}u_{R}\right\vert ^{2}%
+E_{0}u_{R}^{2}\right)  dx}{\int u_{R}^{2}dx}+\frac{\int N(u_{R})dx}{\int
u_{R}^{2}dx}\leq\nonumber\\
&  \leq E_{0}+\frac{N(s_{0})}{s_{0}^{2}}+O\left(  \frac{1}{R}\right)
\label{casto}%
\end{align}

Then by (\ref{W1}) we can easily deduce that for $R$ large enough we have
\begin{equation}
\Lambda(u_{R})<E_{0}. \label{astra}%
\end{equation}

Finally by (\ref{astra}) and Corollary \ref{interm} we get
\[
\Lambda(u_{R})<\Lambda_{0}%
\]

$\square$

\bigskip

Now we are ready to prove the theorem.\bigskip

\textbf{Proof of Theorem \ref{solitoni} } We just need to check that all the
assumptions of Th. \ref{astra1} are satisfied. (EC-1),(EC-2) hold
trivially\textit{. }(EC-3) and (EC-4) hold by Lemma \ref{EC3} and \ref{EC4}
respectively and (\ref{hh}) is verified in Lemma \ref{seguente2}.

$\square$

\section{ Hylomorphic solitons for the generalized BO equation\label{HSK}}

In this section we will study equation (\ref{GBO}).

\begin{proposition}
Let $W$ be a $C^{1}$ function and $u(t,\cdot)\in H^{2s}$ be a smooth solution
of equation (\ref{GBO}) decayng sufficiently fast in space. Then $u$ has the
following constants of motion: the energy
\begin{equation}
E=\int\left(  \frac{1}{2}\left[  D^{s}u\right]  ^{2}+W(u)\right)
dx\label{ebis}%
\end{equation}
and the charge%
\begin{equation}
C=\frac{1}{2}\int u^{2}dx\label{cbis}%
\end{equation}

\end{proposition}

\textbf{Proof.} We have%

\begin{align*}
dE(u)\left[  v\right]   &  =\int\left(  D^{s}u\ D^{s}v+W^{\prime}(u)v\right)
dx\\
&  =\int\left[  D^{2s}u+W^{\prime}(u)\right]  vdx
\end{align*}
hence, using the equation (\ref{GBO})
\begin{align*}
\frac{d}{dt}E(u(t))  &  =\int\left[  D_{x}^{2s}u+W^{\prime}(u)\right]
\partial_{t}u\\
&  =\int\left[  D_{x}^{2s}u+W^{\prime}(u)\right]  \partial_{x}\left[
D_{x}^{2s}u+W^{\prime}(u)\right] \\
&  =\frac{1}{2}\int\partial_{x}\left(  D_{x}^{2s}u+W^{\prime}(u)\right)
^{2}dx=0
\end{align*}

Then $E$ is constant along the solution $u.$

Let us now show that also $C$ is constant along $u.$ By (\ref{GBO}) we have%
\begin{align}
\frac{d}{dt}C(u)  &  =\int u\partial_{t}udx=\int u\partial_{x}\left[
D_{x}^{2s}u+W^{\prime}(u)\right]  dx\label{si}\\
&  =\int u\partial_{x}D_{x}^{2s}udx+\int u\partial_{x}W^{\prime}(u)dx
\end{align}
Let us compute each piece separately:%
\begin{align}
\int u\partial_{x}D_{x}^{2s}udx  &  =\int u\partial_{x}D_{x}^{s}D_{x}%
^{s}udx=\int uD_{x}^{s}\partial_{x}D_{x}^{s}udx\label{su}\\
&  =\int D_{x}^{s}u\partial_{x}D_{x}^{s}u\ dx=\frac{1}{2}\int\partial
_{x}\left(  D_{x}^{s}u\right)  ^{2}=0
\end{align}
Moreover%
\begin{align}%
{\displaystyle\int}
\partial_{x}W^{\prime}(u)udx  &  =-\int W^{\prime}(u)\partial_{x}%
udx\label{so}\\
&  =-\int\partial_{x}W(u)dx=0
\end{align}

Substituting (\ref{so}) and (\ref{su}) in (\ref{si}) we get%
\[
\frac{d}{dt}C(u)=0
\]
$\square$

We will apply the abstract theory of section \ref{sws} taking $X=H^{s}%
(\mathbb{R})$; in this case a function $u(t,\cdot)\in H^{s}$ is a weak
solution of (\ref{GBO}) if $\forall\varphi\in\mathfrak{D}(\mathbb{R})$
\begin{equation}
\int\left[  \partial_{t}u\varphi+u\partial_{x}D_{x}^{2s}\varphi+W^{\prime
}(u)\partial_{x}\varphi\right]  dx=0\ \label{wbo}%
\end{equation}

\begin{theorem}
\label{theoKdV} Let all the assumptions of Theorem (\ref{solitoni}) hold,
$\gamma$ being here the evolution operator for eq. (\ref{GBO}). Then the
equation (\ref{GBO}) admits a family of hylomorphic solitons $u_{\delta}$
$,\ \delta\in\left(  0,\delta_{\infty}\right)  ,u_{\delta}\in H^{2s}.$
\end{theorem}

\textbf{Proof}: The proof of this theorem is essentially the same as the proof
of Th.\ref{solitoni}. The reason for this relies on the fact that the energy
and the charge for eq. (\ref{NSV}) given by (\ref{ghj}) and (\ref{chery}) are
formally the same as the energy and the charge of equation (FKdV) given by
(\ref{ebis}) and (\ref{cbis}). The fact that in the first case $\psi$ is
complex while in the second case $u$ is real-valued does not affect the estimates.

$\square$

The next theorem is the analogous of Th. \ref{marina} and it gives more
information on the structure of the solitons and of its dynamics. Moreover it
permits to eliminate assumption (\ref{Wb}).

\begin{theorem}
\label{marina+}Let all the assumptions of Theorem \ref{marina} hold, $\gamma$
being here the evolution operator for eq. (\ref{GBO}). Then the equation
(FKdV) admits a family of hylomorphic solitons $u_{\delta},\ \delta\in\left(
0,\delta_{\infty}\right)  ,u_{\delta}\in H^{2s}..$ Moreover, $u_{\delta}$ is a
(weak) solution of the equation%
\begin{equation}
\partial_{x}D^{2s}u+\partial_{x}W^{\prime}(u)=\lambda_{\delta}\partial_{x}u
\label{sat+}%
\end{equation}
and
\begin{equation}
U_{\delta}(t,x)=u_{\delta}(x-\lambda_{\delta}t)
\end{equation}
solves (\ref{GBO}).
\end{theorem}

\textbf{Proof}. First let us assume (\ref{Wb}).Then, by Theorem \ref{theoKdV},
the equation (\ref{GBO}) admits a family of hylomorphic solitons $u_{\delta}$
$,\ \delta\in\left(  0,\delta_{\infty}\right)  ,u_{\delta}\in H^{2s}.$

Since $u_{\delta}$ is a minimizer of the energy $E$ on the manifold
$\mathfrak{M}_{c_{\delta}},$ there exists a Lagrange multiplier $\lambda
_{\delta}$ s.t.%

\[
E^{\prime}(u_{\delta})=\lambda_{\delta}C^{\prime}(u_{\delta}).
\]
The above equality can be written as follows%

\[
D^{2s}u+W^{\prime}(u_{\delta})=\lambda_{\delta}u_{\delta}%
\]
So, if we take the derivative $\frac{\partial}{\partial x}$ on both side, we
get (\ref{sat+}). Finally (\ref{sat+}) implies that the travelling wave
$U_{\delta}(t,x)=u_{\delta}(x-\lambda_{\delta}t)$ solves (\ref{GBO}) and
consequently $u_{\delta}$ is a soliton.

It remains to show that the same result holds even when (\ref{Wb}) is
violated. Consider the following equation
\begin{equation}
\partial_{t}u+\partial_{x}\left[  D_{x}^{2s}u+W_{0}^{\prime}(u)\right]  =0,
\label{aaa}%
\end{equation}
where $W_{0}^{\prime\prime}(0)=-2E_{0}<0$. In this case is convenient to
consider the equation%
\begin{equation}
\partial_{t}v+\partial_{x}\left[  D_{x}^{2s}v+W^{\prime}(v)\right]  =0,
\label{bbb}%
\end{equation}
where $W(r)=W_{0}(r)+\left(  1+E_{0}\right)  r^{2}.$We have that
\[
W^{\prime\prime}(0)=2>0
\]
and to every solution $v$ of eq. (\ref{bbb}) corresponds a solution%
\[
u(t,x)=v(t,x-\lambda t)\ \ with\ \ \lambda=2\left(  1+E_{0}\right)
\]
of eq. (\ref{aaa}). In fact%
\begin{align*}
\partial_{t}u+\partial_{x}\left[  D_{x}^{2s}u+W_{0}^{\prime}(u)\right]   &
=\partial_{t}v-\lambda\partial_{x}v+\partial_{x}\left[  D_{x}^{2s}%
v+W_{0}^{\prime}(v)\right] \\
&  =\partial_{t}v-\lambda\partial_{x}v+\partial_{x}\left[  D_{x}%
^{2s}v+W^{\prime}(v)+2\left(  1+E_{0}\right)  v\right] \\
&  =\partial_{t}v-\lambda\partial_{x}v+\partial_{x}\left[  D_{x}%
^{2s}v+W^{\prime}(v)+\lambda v\right] \\
&  =\partial_{t}v-\partial_{x}\left[  D_{x}^{2s}v+W^{\prime}(v)\right]  =0
\end{align*}

$\square$

\bigskip

\begin{remark}
Th. \ref{marta} is a particular case of the above theorem when $D^{2s}$
reduces to a differential operator. Th. \ref{martina} is obtained by
Th.\ref{marina+} taking $s=1/2$ in (\ref{GBO}).
\end{remark}

\begin{remark}
The Cauchy problem for (\ref{BO}) is globally well posed in $H^{1}%
$\cite{Tao2004}. $,$ whereas the Cauchy problem for KdV equation (\ref{1K}) is
globally well posed in $H^{2}$ (see \cite{kato}, \cite{We86}).
\end{remark}

\end{document}